\documentclass[12pt]{article}

\usepackage{amsmath}
\usepackage{multirow} 

\usepackage[cp850]{inputenc}
\usepackage{latexsym,amsfonts}

\usepackage{graphicx}
\usepackage{subcaption}
\usepackage{pdflscape}

\usepackage{caption}

\usepackage{tikz-cd} 
\usepackage{tikz}

\usepackage{enumitem} 
\usepackage{authblk}

\newtheorem{Thm}{Theorem}[section]
\newtheorem{Note}{Note}
\newtheorem{Lem}[Thm]{Lemma}

\newtheorem{Con}[Thm]{Conjecture}

\author[1]{Mustafa Gezek}
\author[2]{Rudi Mathon}
\author[3]{Vladimir D. Tonchev}
\affil[1]{Department of Mathematics, Tekirdag Namik Kemal University, 
Tekirdag, Turkey 59030, mgezek@nku.edu.tr}

\affil[2]{Department of Computer Science, University of Toronto, Toronto,
Ontario, Canada M5S3G4,
rmathon@gmail.com}
\title{\bf Maximal arcs, codes, and new links between projective planes of order 16}
\affil[3]{Department of Mathematical Sciences, Michigan Technological 
University, Houghton, MI USA 49931, tonchev@mtu.edu}

\title{\bf Maximal arcs, codes, and new links between projective planes}

\begin{document}
\maketitle

\begin{abstract}

In this paper we consider binary linear codes spanned by incidence matrices of
Steiner 2-designs associated with maximal arcs in projective planes
of even order, and their dual codes.
Upper and lower bounds on the 2-rank of the incidence matrices
are derived.
A lower bound on the minimum distance of the dual codes is proved, and
it is shown that the bound is achieved if and only if the
related maximal arc contains a hyperoval of the plane.
The  binary linear codes
of length 52 spanned by the incidence matrices of 2-$(52,4,1)$ designs
associated with previously known and some newly found 
maximal arcs of degree 4 in projective planes of order 16
are analyzed and classified up to equivalence.
The classification shows that some designs associated with maximal
arcs in nonisomorphic planes generate equivalent codes.
This phenomenon establishes new links between several of the
known planes. A conjecture concerning the codes
of maximal arcs in $PG(2,2^m)$ is formulated.
 
\medskip


\end{abstract}

\vspace{5mm}

\section{Introduction}

We assume familiarity with the basic facts and notions 
from design theory, finite geometries, and coding theory 
\cite{BJL,BBFKKW,CRC,Ding,Hir}.

A 2-$(v,k,\lambda)$ design (or shortly, a 2-design) is a pair $D$=$\{ X, B \}$
of a set $X$ of $v$ {\it points} and a collection $B$ of subsets
of $X$ of size $k$ called {\it blocks}, such that every two points
 appear together in
exactly $\lambda$ blocks. 
Every point of a
$2$-$(v,k,\lambda)$ design is contained in
$r=\lambda(v-1)/(k-1)$ blocks, and the total number
of blocks is $b=v(v-1)\lambda/k(k-1)$.

The incidence matrix of a design $D$ is
 a $(0,1)$-matrix $A=(a_{ij})$ with rows
labeled by the blocks, columns labeled by the points, where
$a_{i,j}=1$ if the $i$th block contains the $j$th point, and $a_{i,j}=0$
otherwise. If $p$ is a prime, the $p$-rank of a design $D$ is the rank
of its incidence matrix over a finite field of characteristic $p$.
Two designs are {\it isomorphic} if there is a bijection between their
point sets that maps every block of the first design to a block of the second
design. An {\it automorphism} of a design is any isomorphism of the design
to itself.
The set of all automorphisms of $D$ form the automorphism group
$Aut(D)$ of $D$.
The {\it dual} design $D^\perp$ of a design $D$
has as points the blocks of $D$, and as blocks the points of $D$.
A 2-$(v,k,\lambda)$ design is {\it symmetric} if $b=v$,
or equivalently,  $r=k$.
The dual design $D^\perp$ of a symmetric 2-$(v,k,\lambda)$ design $D$
is a symmetric design with the same parameters as $D$.
A symmetric design $D$ is {\it self-dual} if $D$ and $D^{\perp}$ are
isomorphic.

A design with $\lambda=1$ is called a {\it Steiner} design.
An affine plane of order $n$ ($n\ge 2$),
 is a Steiner 2-$(n^2, n, 1)$ design.
A projective plane of order $n$ is a symmetric
Steiner 2-$(n^2 + n +1, n+1, 1)$ design with $n\ge 2$.
The classical (or Desarguesian) plane $PG(2,p^t)$
of order $n=p^t$, where $p$ is prime and $t\ge 1$, has as points
the 1-dimensional subspaces of the 3-dimensional vector space $V_3$
over the finite field of order $p^t$, and as blocks (or {\it lines}), the
2-dimensional subspaces of $V_3$.

Let $D=\{ X, B \}$ be a Steiner 2-$(v,k,1)$ design with
point set $X$, collection of blocks $B$, and let $v$ be a multiple
of $k$,  $v=nk$. Since every point of $X$ is contained in
$r=(v-1)/(k-1)=(nk-1)/(k-1)$ blocks,  $k-1$ divides $n-1$.
Thus, $n-1=s(k-1)$ for some integer $s\ge 1$, and
\[ v=nk = (sk -s +1)k. \]
A {\it parallel class} of $D$ is a set of $v/k=n$
pairwise disjoint blocks,
and a {\it resolution} of $D$ is a partition of the collection of blocks
$B$ into $r=(v-1)/(k-1)=sk+1$ disjoint parallel classes.
A design is {\it resolvable} if it admits a resolution.

Any 2-$((sk -s +1)k, k, 1)$ design with $s=1$ is
an affine plane of order $k$, and admits exactly one resolution.
If $s>1$, a resolvable 2-$((sk -s +1)k, k, 1)$ design may admit more than one resolution.

Let $P$ be a projective plane of order $q$, and let $m$ and $k$ be positive integers such that $k\le q+1$ and $k \le m \le q^2 + q + 1$. An $(m,k)$-{\it arc} (or an arc of size $m$ and degree $k$) is a set $A$ of $m$ points such that every line of $P$ contains at most $k$ points from $A$.
 An arc of degree 2 is also called an {\em oval}. Let $x$ be a point in an 
 $(m,k)$-arc $A$, let $L_1,\ldots, L_{q+1}$ be the lines through $x$, 
 and let $n_i = |A\cap L_i|$, $1\le i \le q+1$.
Then\[  m = 1 +  \sum_{i=1}^{q+1} (n_i -1). \] Since $n_i \le k$, it follows that
\[ m \le qk+k -q, \] and the equality $m = qk+k -q$ holds if and only if every line of $P$ is either disjoint from $A$ or meets $A$ in exactly $k$ points.
An $(m,k)$-arc  is called {\it maximal} if $m = qk+k -q$. 
A {\it hyperoval} is a maximal arc of degree 2.

Let $A$ be a maximal $(qk+k-q,k)$-arc with $k\le q$ ,
 thus $m=qk+k-q \le q^2$, and let $y$ be a point outside $A$. 
 Let $n$ be the number of lines that meet 
$A$ in $k$ points. We have
\[ nk = qk + k -q. \] 
Thus, $k$ divides $q$, and $q=sk$ for some integer $s\ge 1$
\footnote{In all known exmples $q$ is a poewr of 2.}.
If $q=sk$, the size of a maximal arc of degree $k$ can be written
in terms of $s$ and $k$ as $(sk-s+1)k$.

The set of lines that are disjoint from a maximal $((sk-s+1)k,k)$-arc 
 $A$ in a projective plane $P$ of order $q=sk$, form a maximal
$((sk-k+1)s,s)$-arc $A^\perp$ in the dual plane $P^\perp$,
called the {\it dual arc} of $A$.

Maximal arcs with $1<k < q$ do 
not exist in any Desarguesian plane of odd order
$q$ \cite{BBM},
and are known to exist in every Desarguesian plane of
order $q=2^t$ \cite{Den}, \cite{HM1}, \cite{HM2}, 
\cite{Math}, \cite{Thas74},
 as well as in some non-Desarguesian planes of even order
 $q=2^t$ 
\cite{Gezek}, \cite{GTW}, \cite{H1}, \cite{H2}, \cite{H3}, 
\cite{HST}, \cite{PRS}, \cite{Thas80}.

If $k>1$, the non-empty intersections of
a maximal $((sk-s+1)k, k)$-arc $A$ with 
lines of a projective plane $P$ of order $q=sk$
are the blocks of a resolvable 2-$((sk-s+1)k,k,1)$ design $D$.
Similarly, if $s>1$, the dual $((sk-k+1)s, s)$-arc ${A}^{\perp}$
in the dual plane $P^\perp$
is the point set of a resolvable
2-$((sk-k+1)s, s, 1)$ design $D^{\perp}$.
We will refer to $D$ (resp. $D^{\perp}$) as a design embeddable in
$P$ (resp. ${P}^\perp$) as a maximal arc.

Two maximal arcs $A'$, $A''$ in a projective plane
$P$ are {\it equivalent }
if there is a collineation of $P$ that maps $A'$ to  $A''$.
Designs associated with equivalent arcs
are necessarily isomorphic, while the converse is not true in general.

Let $D$ be a resolvable Steiner 2-$(v,k,1)$ design.
Two resolutions $R_1,  R_2$ of $D$,
where
\begin{equation}
\label{def}
 R_1 = P^{(1)}_1 \cup P^{(1)}_2 \cup \cdots P^{(1)}_r, \  R_2 =
 P^{(2)}_1 \cup P^{(2)}_2 \cup \cdots P^{(2)}_r,
\end{equation}
are called
{\em compatible} \cite{ton-res}, if they share one parallel class,
$ P^{(1)}_i = P^{(2)}_j$,
and
\[ |P^{(1)}_{i'} \cap P^{(2)}_{j'}|\le 1 \]
for  $(i',j') \neq (i,j)$.

The following theorem
gives an upper bound
on the number of pairwise compatible resolutions of a
resolvable 2-$((sk-s+1)k, k, 1)$ design,
and characterizes the designs for which this upper bound is achieved.
\begin{Thm}
\cite{ton-res}.
\label{th}
Let $S=\{ R_1,\ldots, R_m \}$ be a set of $m$ mutually  compatible resolutions
of a  2-$((sk -s +1)k, k, 1)$ design $D=\{ X, B \}$.
Then
\[ m\le (sk-k+1)s. \]
The equality
\[ m=(sk-k+1)s \]
holds if and only if there exists a projective plane ${P}$ of order $q=sk$
such that $D$ is embeddable in ${P}$ as a maximal $((sk-s+1)k, k)$-arc.
\end{Thm}

In Section \ref{secN2},
we consider binary linear codes spanned by
the rows of incidence matrices of
Steiner 2-designs associated with maximal arcs in projective planes
of even order, and their dual codes.
Upper and lower bounds on the 2-rank of the incidence matrices
are derived.
A lower bound on the minimum distance of the dual codes is proved, and
it is shown that the bound is achieved if and only if the
related maximal arc contains a hyperoval of the plane.

In Section \ref{secn}, we analyze the binary linear codes
of length 52 spanned by the incidence matrices of 2-$(52,4,1)$ designs 
associated with maximal arcs of degree 4 in projective planes of order 16.
The codes associated with maximal arcs in $PG(2,16)$ are
distance optimal, while one code associated with an arc in the semi-filed
plane is shown to be optimal with respect to
each of its parameters: minimum  distance, dimension, and length. A conjecture concerning the codes
of maximal arcs in $PG(2,2^m)$ is formulated.

The codes are classified according to their dimension, 
and all codes having the same dimension are further classified up to
equivalence.
The classification shows that some designs associated with maximal
arcs in nonisomorphic planes generate equivalent codes.
This phenomenon establishes new links between several of the
known planes that are discussed in Section \ref{sec3}.

Section \ref{sub2} lists explicitly eleven new maximal arcs
of degree 4 found recently in the planes DEMP, SEMI2, LMRH, HALL, and BBH1.

\section{Binary codes of designs arising from maximal arcs}
\label{secN2}

An {\it oval} of a Steiner 2-$(v,k,1)$ design $D$ with $k \ge 2$
 is a set $S$
that meets every block in at most two points.

\begin{Lem}
\label{lem1}
(a) The size of an oval $S$ of a Steiner 2-$(v,k,1)$ design
is bounded above by
\[ |S| \le r+1, \]
where $r=(v-1)/(k-1)$ is the number of blocks through a point.\\
(b) The equality $|S|=r+1$ holds if and only if every block
is either disjoint from $S$, or meets $S$ in exactly two points.
\end{Lem}
\begin{Proof}
Let $n_i$ ($i=0, 1, 2$) denote the number of blocks meeting $S$
in $i$ points. Counting in two ways the incident pairs of points
from $S$ and blocks, and incidence pairs of pairs of points
from $S$ and blocks, we have
\begin{eqnarray*}
 n_1   +   2n_2 &  = & |S|r, \\
        2n_2 & = & |S|(|S|-1),
\end{eqnarray*}
hence
\[ n_1 = |S|r - |S|(|S|-1) \ge 0, \]
and the statement follows. {$\Box$}
\end{Proof}

An oval of size $r+1$ is called
 a {\it hyperoval}.

\begin{Note}
\label{not1}
{\rm
If $D$
is a symmetric Steiner design, that is, a projective plane,
a hyperoval is simply a maximal arc of degree two. 
}
\end{Note}

\begin{Thm}
\label{t1.2}
Let $C$ be a binary linear code of length $n=2^{m+s}-2^m + 2^s$,
 spanned by
the rows of the incidence matrix $A$
of a Steiner 2-$(2^{m+s}-2^m + 2^s, 2^s, 1)$ design $D$
with $m\ge s\ge 1$. \\
(i) The all-one vector $\bar{1}=(1,\ldots,1)$ belongs to $C\cap C^\perp$.\\
(ii) The dual code $C^\perp$ admits majority-logic decoding that corrects
up to $t=2^{m-1}$ errors.\\
(iii) The minimum distance $d^\perp$ of $C^\perp$
is an even number equal to $2^m +2$ if $D$ contains 
hyperovals, and
$d^\perp \ge 2^m + 4$ if $D$ has no hyperovals.\\
(iv) The minimum distance $d$ of $C$ is an even number smaller than or equal to $2^s$.\\
(v) The dimension $k$ of $C$, or equivalently, the
2-rank of $A$, $rank_{2}A$, satisfies the inequalities
\[ 1+ \lceil \log_{2}(\sum_{i=0}^t {n-1 \choose i} )\rceil \le 
rank_{2}A \le 
n - 1 - \lfloor \log_{2}(\sum_{i=0}^{d/2 -1} { n-1 \choose i})\rfloor,\]
where $t= 2^{m-1}$ if $d^\perp = 2^m +2$,
 and $t = d^{\perp}/2 -1$ if $d^\perp \ge 2^m + 4$.

\end{Thm}
\begin{Proof}
(i) Since all rows of $A$ are of even weight $2^s$, every
row is orthogonal to $\bar{1}$ over $GF(2)$, hence
$\bar{1}\in C^\perp$.
 Every point of $D$ is contained in $r=2^m +1$ blocks.
Thus, every column of $A$ contains $2^m +1$ nonzero entries,
and the binary sum of all rows of $A$ is equal to  $\bar{1}$,
hence $\bar{1}\in C$.\\
(ii) The rows of $A$ provide a set of checks that can be used
to correct up to
\[ t = \lfloor{\frac{r}{2}}\rfloor= \lfloor{\frac{2^m +1}{2}}\rfloor = 
2^{m-1} \]
errors in $C^\perp$ by using the Rudolph majority-logic decoding
algorithm (cf. \cite{Ru},
\cite[Theorem 8.1, page 1252]{Thbc}).\\
(iii) It follows from (ii) that $d^\perp \ge 2t+1=2^m + 1$.
Assume that $d^\perp = 2^m + 2$, and let $S$ be the support
of a minimum weight codeword in $C^\perp$. Clearly, $S$ meets
every block of $D$ in an even number of points. Let $n_{2i}$
denote the number of blocks that meet $S$ in $2i$ points,
$0 \le i \le 2^{s-1}$. Counting in two ways the occurrences of single points, and ordered pairs of points of $S$ in blocks 
of $D$, we have
\begin{eqnarray*}
\sum_{i=1}^{2^{s-1}} 2in_{2i} & = & (2^m + 2)(2^m +1),\\
\sum_{i=1}^{2^{s-1}} 2i(2i -1)n_{2i} & = & (2^m + 2)(2^m +1).
\end{eqnarray*}

Subtracting the first of the above equations from the second gives
\[ \sum_{i=2}^{2^{s-1}} 4i(i-1)n_{2i} =0, \]
hence $n_{2i}=0$ for $i>1$, and $S$ is a hyperoval of $D$.
Consequently, by Lemma \ref{lem1},
the number of codewords of $C^\perp$
having weight $2^m +2$ is equal to the number of hyperovals of $D$.
If $D$ has no hyperovals, then $d^\perp > 2^m + 2$, and since all
weights in $C^\perp$ are even by part (i), it follows that
$d^\perp \ge 2^m + 4$.

Part (iv) is obvious.
The upper bound in Part (v) follows from
applying the sphere packing bound (cf., e.g. \cite[1.12]{HP}),
 to a punctured $[n-1,k,d -1]$ code $C'$ of $C$.
The lower bound in Part (v) follows from applying the sphere packing bound
to a punctured $[n-1,n-k,d^\perp -1]$ code of $C^\perp$.  
$\Box$

\end{Proof}

\begin{Note}
\label{n2}
{\rm
According to Theorem \ref{t1.2} (ii) and (iii), if $D$ contains 
hyperovals
then the dual code $C^\perp$ can correct
the maximum number of errors guaranteed by its minimum distance,
by using efficient majority-logic decoding,
and in addition, the number of codewords in $C^\perp$ having  minimum weight 
is equal to the number of hyperovals of $D$.
}
\end{Note}

As a corollary of Theorem \ref{t1.2}, we have the following.

\begin{Thm}
\label{t1.3}
Let $C$ be a binary code of length $n=2^{m+s}-2^m + 2^s$,
 spanned by
the incidence matrix $A$
of a Steiner 2-$(2^{m+s}-2^m + 2^s, 2^s, 1)$ design $D$
with $m\ge s\ge 1$,
associated with a maximal arc $S$ of degree $2^s$ in a projective
plane $\Pi$ of even order $q=2^m$.
Let $d$ and $d^\perp$ denote the minimum distance of $C$ and $C^\perp$
respectively.\\

(a) The minimum distance $d^\perp$ of the dual code $C^\perp$ is $d^\perp = 2^m + 2$ if
$D$ contains a hyperoval of $\Pi$, and $d^\perp \ge 2^m + 4$ otherwise.\\

(b) The minimum distance $d$ of $C$ is an even number smaller
than or equal to $2^s$.
 If $d=2^s$ then
\[ 1 + \lceil  \log_{2}(\sum_{i=0}^t {n-1 \choose i} )\rceil \le
rank_{2}A \le n - 1 - \lfloor \log_{2}(\sum_{i=0}^{2^{s-1} -1}
{ n-1 \choose i})\rfloor, \]
where where $t= 2^{m-1}$ if $D$ contains a hyperoval of
$\Pi$, and  $t = d^{\perp}/2 -1 \ge 2^{m-1} +1$ if $D$ does not contain any hyperoval
of $\Pi$.

\end{Thm}

\begin{Proof}
 Since every hyperoval of $D$  is also a hyperoval of $\Pi$,
(a) and (b) follow from part (iii) and (v) of Theorem \ref{t1.2}
respectively.
$\Box$
\end{Proof}

We note that if $s=1$, the design $D$ from Theorems \ref{t1.2} and \ref{t1.3}
is the trivial design on $2^m + 2$ points having as blocks
all pairs of points. In this case, $rank_{2}A= 2^m +1$, and the code
$C$ is the unique binary code consisting of all  vectors
of even weight, while $C^{\perp} = C\cap C^\perp  =\{ {\bf 0}, \bar{1} \}$,
where ${\bf 0}$ stands for the zero vector.
Apart from this trivial case, the question about the 2-rank of designs arising
from maximal arcs appears to be widely open, with one notable exception concerning
the case $s=m-1$: it was proved by Carpenter \cite{Ca} that the 2-rank of a
2-$(2^{2m-1} - 2^{m-1}, 2^{m-1}, 1)$ design associated with the dual arc
of a regular hyperoval in $PG(2,2^m)$ is equal to $3^m - 2^m$.

\section{Designs associated with maximal arcs of degree 4    
and their codes}
\label{secn}

A Steiner design $D$ associated with a maximal arc of degree 4 
in a plane of order $2^m$ has parameters
2-$(3\cdot 2^m + 4, 4, 1)$, and its binary code $C$ is of length 
$n=3\cdot 2^m + 4$
and has minimum distance $d=2$ or $d=4$.
In this section,
we consider 2-$(3\cdot 2^m + 4, 4, 1)$ designs and their 
binary codes
when  $m=2$, 3 and 4.

If  $m=2$, $D$ is the unique (up to isomorphism)
2-$(16,4,1)$ design, being isomorphic to the affine plane of order 4,
and can be viewed  as the design associated with a maximal arc of
degree 4 in the projective plane of order 4.
The 2-rank of its incidence matrix $A$ is 9, hence the code $C$
spanned by $A$ has dimension 9.
The minimum distance of $C$ is 4, and there are exactly 20 codewords
of minimum weight 4, being the rows of the incidence matrix 
$A$.
The design $D$ contains ovals, hence $d^\perp = 6$ by Theorem 1.2, (iii).
The lower bound and upper bound on $rank_{2}A$ from Theorem \ref{t1.2}
are 8 and 11 respectively.
The code $C$ is distance optimal in the sense that
4 is the largest possible minimum distance for a binary
linear code of  length 16  and dimension 9.
The dual $[16,7,6]$ code $C^\perp$ is also distance optimal
(see \cite{Gr}).

In the next case, $m=3$, the design parameters are 
2-$(28,4,1)$.
There are at least 4653 known nonisomorphic designs with these parameters
\cite{BBT}, \cite{Kr}, having 2-ranks ranging from 19 to 26 \cite{BBT},
with the exception of 2-rank 20 (it was shown in \cite{JaffeT}
that there are no 2-$(28,4,1)$ with 2-rank 20).
It was proved in \cite{MTW} that the minimum 2-rank
of any 2-$(28,4,1)$ design is 19, and up to isomorphism, there is a
unique design of minimum 2-rank 19, being isomorphic to the design $D$
associated with a maximal $(28,4)$-arc in the projective plane of order 8,
$PG(2,8)$ (also referred to as the Ree unital \cite{Br}).
The minimum distance of the code $C$
of the Ree unital $D$ is 4, and every codeword of minimum weight
corresponds to a block of $D$. The Ree unital,
being an oval design in the terminology of \cite{Ca}, 
contains ovals that are hyperovals
in the projective plane of order 8, hence  the minimum distance of
the dual code $C^\perp$ is $d^\perp = 10$, and $C^\perp$ can correct
up to $\lfloor (10 -1)/2 \rfloor = 4$ errors by majority-logic decoding.
The $[28,19,4]$ code $C$ 
and its dual $[28,9,10]$ code $C^\perp$ are both distance optimal
(cf. \cite{Gr}).


If $m=4$, the parameters of the design $D$ from 
Theorems \ref{t1.2} and \ref{t1.3}
 are 2-$(52,4,1)$ and correspond to a design 
associated with a maximal $(52,4)$-arc in a projective plane of order 16.

Up to isomorphism, there are 22 known projective planes of order 16.
The only projective plane of order 16 for which
all inequivalent maximal arcs of degree 4 have been
completely classified is the Desarguesian plane $PG(2,16)$.
There are exactly two inequivalent maximal $(52,4)$-arcs in $PG(2,16)$
\cite{BB}, with collineation groups of order 68 and 408.
The arc with group of order 408 admits a cyclic collineation of
order 51 \cite{DDT}.
The maximal arcs of degree 4 have not been classified completely
in any of the the 21 known  non-Desarguesian planes
of order 16.  Maximal arcs of degree 4 have been found in
all but four of the known non-Desarguesian planes of order 16
\cite{Gezek}, \cite{GTW}, \cite{HST}.
All previously known maximal $(52,4)$-arcs 
can be found in \cite{GTW} and \cite{HST}.
 Eleven new maximal arcs
of degree 4 are given in Section \ref{sub2} of this paper.
The line sets of the projective planes of order 16
and all known maximal arcs of degree 4,
including the new arcs desribed in Section \ref{sub2},
 are available online
at
http://pages.mtu.edu/$\sim$tonchev/planesOForder16.txt
and
http://pages.mtu.edu/$\sim$tonchev/pointsetsOFmaxArcs.txt
respectively.

Table \ref{tab1} contains data about the 2-$(52,4,1)$ designs associated 
with previously known and the eleven newly found
maximal arcs of degree 4
described in Section \ref{sub2}.
 Column 3 lists the orders of the stabilizers of the maximal arcs,
which happen to be also the full automorphism groups of the related designs.
 Column 4 gives the number of parallel classes of the design $D$
associated with the given arc, followed by the number of parallel
classes of the design $D^{\perp}$ associated with the dual arc.
Column 5 lists the number of resolutions of $D$ and $D^{\perp}$.
The last column lists the orders of the automorphism groups
of the codes $C(D)$ and $C(D^{\perp})$ respectively.
Cliquer \cite{cliquer} and Magma \cite{magma} were
used for these computations.

Since the number of parallel classes, resolutions, code dimension
and minimum distance are invariant under isomorphismis of designs,
the data from Table \ref{tab1} implies that the number of pairwise
nonisomorphic resolvable 2-$(52,4,1)$ designs associated with maximal arcs
is greater than or equal to 55. Further computation
of possible isomorphisms using Magma \cite{magma} shows that the number of 
nonisomorphic designs is exactly  55.

We note that the previously known lower bound on the number of nonisomorphic
resolvable 2-$(52,4,1)$ designs given in \cite{CRC} is 30.

 The 2-ranks of these 55 designs vary from 41 to 49, and the minimum 41
is achieved only by designs associated with maximal arcs in the Desarguesian 
plane $PG(2,16)$.

Table \ref{tab3} lists the 2-ranks of the designs,
and gives a partition of the codes into equivalence classes,
the numbers $A_2$ and $A_4$ of
codewords of weight 2 and 4 respectively,
and the orders of the automorphism groups of the codes.
Table \ref{tab4} lists the parameters of the codes and their dual codes,
as well as the number $A_{d^\perp}$ of codewords of minimum weight in
$C^\perp$. We note that according to Theorem \ref{t1.3} $A_{d^\perp}$
is equal to the number of hyperovals of the plane contained in the point set of
$D$ in the cases when $d^\perp =18$.

Four of the 55 codes associated with these designs have minimum
distance $d=4$: the codes $C'$, $C''$ of the two 
maximal arcs in $PG(2,16)$, and the codes
$C'''$ and $C^{iv}$  of dimension 45 associated with the maximal arcs
 SEMI4.1 and SEMI2.7 (see Table \ref{tab1} and Table \ref{tab4}).
Every codeword of minimum weight in $C'$ and $C''$
is the incidence vector of a block, while the number $A_4 = 4469$  of
minimum weight codewords in $C'''$ and $C^{iv}$
is larger than the number of blocks.
These observations suggest that the following statement may be true
in general.
\begin{Con}
\label{conj2}
The binary code spanned by the incidence matrix $A$
of a design $D$ associated with a maximal arc of degree 4 in $PG(2,2^m)$,
where $m\ge 2$, has minimum distance 4, and every codeword of
 minimum weight is a row of $A$.
\end{Con}

 We note that the codes $C'''$ and $C^{iv}$
are equivalent (see Table \ref{tab3} that lists the equivalence 
classes of codes),
while the codes $C'$, $C''$ of the two arcs in $PG(2,16)$ are inequivalent, 
despite
the fact that both codes have identical weight distributions.
It follows from \cite{DDT} that one of these codes, namely the code with full 
automorphism group of order 408, is an extended cyclic code.
Both $(52,4)$-arcs in $PG(2,16)$ 
contain hyperovals, hence $d(C')^\perp = d(C'')^\perp = 18$
by Theorem \ref{t1.3}. In addition, the number of codewords
of minimum weight $A_{d^\perp}$ in each of the codes $(C')^\perp$
and $(C'')^\perp$ equals the number of hyperovals contained in the corresponding
maximal arcs, namely $A_{d^\perp} = 54$ (see Table \ref{tab4}).

A comparison with the parameters of best known error-correcting codes \cite{Gr}
shows that the highest known minimum distance of a binary code of length 
52 and dimension
41 is $d=4$, while the best known theoretical upper bound is $d\le 5$.
This makes the $[52,41,4]$ codes associated with the maximal arcs in $PG(2,16)$
"nearly" optimal.

In comparison, the $[52,45,4]$ code
$C'''$ of the arc SEMI4.1 is distance optimal. The dimension 45 of $C'''$
is just by one shorter from the upper bound 46 obtained from Theorem \ref{t1.3}.
However, the minimum distance of any binary $[52,46]$ code is at most 3 
(which is seen by applying the sphere packing bound on a punctured $[51,46,3]$ 
code),
therefore the $[52,45,4]$ code $C'''$ is both distance and dimension optimal.
This implies that if the binary code of a 2-$(52,4,1)$ design $D$ has minimum 
distance 4, the 2-rank of $D$ cannot exceed 45.
In addition, since a binary $[51,45,4]$ code does not exist, 
because the parameters of a punctured $[50,45,3]$ code
violate the sphere packing bound, the $[52,45,4]$ code $C'''$ is also length
optimal, that is, 52 is the smallest possible length for a binary code 
of dimension 45 and minimum distance 4.

\section{New connections between projective planes}
\label{sec3}

Table \ref{tab2} lists all previously known connections,
as well as some new connections described in this section,
between nonisomorphic projective planes of order 16.
An entry in a row and a column labeled by the same
(non self-dual) plane
indicates a connection between the plane and its dual plane that is
based on designs associated with maximal arcs or their codes.

An entry ${\bf 1}$ indicates that the corresponding
planes are connected by derivation \cite{HST},\cite{M}, ${\bf 2}$
indicates that the corresponding
planes are connected by superderivation \cite{MathonTalk},
and ${\bf 3}$ indicates
that the planes share a semibiplane \cite{M}.

An entry ${\bf 4}$ indicates that the planes share a
2-$(52,4,1)$ design associated with a maximal $(52,4)$-arc
via a construction based on Theorem \ref{th}.
Connections of this type were considered in \cite{GTW}.

An entry ${\bf 5}$ indicates a new connection between a pair of
nonisomorphic planes that
share a binary linear code of length 52 generated by
non-isomorphic designs associated
with maximal arcs in the corresponding planes.
Let $D$ be a 2-$(52,4,1)$ design associated with
a maximal arc of degree 4, and let $C(D)$ (resp. $C(D^\perp)$)
 be the binary linear code
of length 52 generated by the incidence matrix of $D$ (resp. $D^\perp$).
The parameters of these codes
are listed in Table \ref{tab4}, 
 and the orders
of their automorphism groups are listed in the last two columns of
Table \ref{tab1}.

The codes were further sorted according to their weight distributions,
and codes having the same weight distribution were tested for equivalences
using Magma \cite{magma}.
This classification shows that the 55 codes
of the 55 nonisomorphic 2-$(52,4,1)$ designs are partitioned
into 27 equivalence classes, listed in Column 3 of Table \ref{tab3}.

Specific permutations that provide equivalences
between codes from the same equivalence
class are listed in Table \ref{tab5}.

The main result implied by this classification of the codes up to equivalence
 is the surprising fact that in several instances the codes of designs 
arising from maximal arcs in  different planes are equivalent, 
hence these codes provide new connections between the corresponding planes
(see lines 4, 5, 6, 11, 13, 17, 18, 20, 24, 26 and 27 in Table \ref{tab2}).
For example, the Mathon plane MATH is linked to the Johnson plane JOHN,
the Lorimer-Rahilly plane LMRH, the    semifield plane  SEMI2, and the
Hall plane HALL.

\section{New maximal $(52, 4)$-arcs }
\label{sub2}

The specific line sets of the known projective planes of order 16
that we are using in this paper were graciously provided to the third author 
by Gordon F. Royle, and are available online at 
\begin{verbatim}
http://pages.mtu.edu/~tonchev/planesOForder16.txt
\end{verbatim}
All previously known maximal arcs of degree 4 are given in 
\cite{GTW} and \cite{HST}, and are available
online at 
\begin{verbatim}
http://pages.mtu.edu/~tonchev/pointsetsOFmaxArcs.txt
\end{verbatim}

Our notation in this section follows \cite{GTW}.
Recently, the
 first author found eight new maximal arcs in 
some of the projective planes of order 16.
The new arcs are unions of orbits of appropriate subgroups of the
the automorphism group of the associated plane.
Two new arcs, stabilized by a nonabelian group of order 12,
were found in the plane DEMP,
and five arcs in the plane SEMI2 and one arc in the plane LMRH were
found as unions of orbits under subgroups of order 16.
The point sets of these eight new maximal arcs are listed below.

\vspace{.1in}
DEMP.3 = $\{263, 265, 266, 258, 32, 122, 142, 243, 187, 102, 61, 197, 84, 232, \\ 210, 156, 18, 126, 140, 251, 181, 112, 52, 195, 88, 237, 214, 154, 30, 117, 144, 244, 
\\ 178, 109, 54, 202, 83, 236, 219, 152, 24, 116, 139, 252, 189, 99, 62, 208, 82, 229, \\ 218, 150\},$

\vspace{.1in}

DEMP.4 = $\{ 
273, 260, 257, 258, 14, 69, 61, 27, 34, 128, 255, 232, 153, 97, 84, \\ 186, 7, 71, 60, 22, 39, 124, 246, 227, 147, 108, 86, 179, 8, 68, 63, 32, 45, 122, 242, \\ 225, 149, 110, 89, 187, 133, 
194, 224, 175, 222, 161, 212, 138, 200, 141, 203, 169
\},$

\vspace{.1in}

SEMI2.3 = $\{ 
263, 268, 265, 267, 23, 28, 228, 25, 27, 234, 240, 229, 4, 16, 124,  \\ 5, 123, 10, 121, 119, 49, 251, 76, 63, 56, 145, 252, 247, 75, 73, 50, 159, 152, 249, 71, \\ 146, 36, 88, 202, 48, 37, 
216, 82, 81, 197, 208, 42, 210, 209, 95, 196, 223
\},$

\vspace{.1in}

SEMI2.4 = $\{ 
259, 269, 262, 270, 20, 233, 32, 21, 231, 236, 26, 235, 4, 16, 124, \\ 5, 123, 10, 121, 119, 33, 90, 208, 47, 40, 209, 85, 96, 196, 202, 34, 223, 216, 84, 197, \\ 210, 49, 251, 66, 63, 56, 
155, 252, 247, 72, 79, 50, 156, 151, 249, 65, 153
\},$

\vspace{.1in}

SEMI2.5 = $\{ 
260, 272, 266, 261, 23, 27, 121, 25, 124, 28, 119, 123, 1, 2, 232, \\ 15, 239, 8, 225, 226, 66, 146, 69, 154, 72, 152, 74, 149, 49, 241, 58, 245, 63, 255, 53, \\ 250, 36, 96, 44, 219, 48, 
89, 199, 84, 224, 43, 212, 220, 196, 87, 208, 201
\},$

\vspace{.1in}

SEMI2.6 = $\{ 
260, 268, 266, 263, 18, 136, 216, 26, 25, 153, 133, 135, 30, 213, \\ 215, 131, 158, 146, 211, 154, 35, 79, 120, 39, 37, 112, 80, 76, 40, 117, 119, 70, 111, \\ 102, 115, 108, 6, 175, 206, 12,
16, 57, 176, 172, 15, 201, 202, 166, 62, 50, 194, 58
\},$

\vspace{.1in}

SEMI2.7 = $\{ 
261, 263, 271, 262, 25, 58, 250, 30, 32, 31, 128, 50, 60, 242, 252, \\ 54, 246, 127, 121, 126, 85, 110, 139, 88, 91, 93, 206, 105, 111, 141, 133, 112, 136, \\ 201, 207, 208, 5, 149, 70, 8, 
11, 13, 226, 152, 155, 76, 66, 157, 74, 234, 230, 236
\},$
\vspace{.1in}

LMRH.2 = $\{ 
46, 78, 250, 90, 42, 74, 94, 254, 260, 266, 270, 269, 20, 29, 132, \\ 141, 4, 13, 164, 173, 25, 27, 50, 194, 137, 145, 209, 139, 9, 11, 146, 210, 169, 49, \\ 193, 171, 37, 70, 64, 195, 69, 
147, 224, 38, 246, 208, 86, 51, 85, 211, 245, 160
\}.$

\vspace{.1in}

A probabilistic search 
algorithm developed by the second author was used to find
three new maximal arcs in the planes DEMP, HALL and BBH1.
Each of these new arcs has a stabilizer of order 4, 
so it is computationally unfeasible to find these arcs with the 
previous method that works well for stabilizers of order at least 12.

A randomised local search was performed to find sets in projective
 planes with prescribed line intersections. At the start
 of a new experiment, a set of points is selected at random.
At each move, a neighbourhood of the current subset, defined
 by a single interchange between a point in the set and a point
 outside, is examined in order to generate a list of swaps
 that minimise an objective function formed from the actual
 and desired line intersections, respectively. A locally optimal
 move is then randomly selected from that list and the set
 updated to reflect it. A tabu list of the most recent moves is
 maintained to prevent cycling, and occasional random moves are
 performed to improve the search efificiency.
 This search algorithm worked well for projective planes of order
 16 and reproduced all the known (52,4)-arcs with frequencies
 inversely proportional to their group orders. We have run typically
 $10^4$ experiments with $5 \times 10^4$ moves per experiment on a MacBook Pro.
 The same approach was used to verify the existence of 4-arcs $PG(2,32)$
 in \cite{Math}. 

 The point sets of these three new maximal arcs are:

\vspace{.1in}

DEMP.5 = $\{ 
1, 3, 8, 15, 23, 24, 25, 28, 36, 38, 41, 43, 51, 54, 61, 64, 66, 69, 70, \\ 78, 81, 82, 89, 96, 100, 104, 106, 109, 129, 133, 139, 140, 149, 154, 156, 160, 178, \\ 183, 189, 190, 195, 202, 
206, 207, 228, 231, 235, 239, 257, 260, 271, 272
\},$

\vspace{.1in}

HALL.2 = $\{ 
1, 2, 5, 14, 19, 27, 28, 32, 34, 39, 40, 45, 49, 53, 54, 63, 81, 84, 88, \\ 95, 103, 107, 108, 109, 131, 134, 142, 143, 147, 153, 154, 155, 166, 169, 170, 173, \\ 180, 183, 184, 185, 197, 
202, 204, 208, 210, 212, 222, 224, 257, 260, 262, 266
\},$

\vspace{.1in}

BBH1.3 = $\{ 
11, 13, 14, 16, 18, 27, 30, 31, 34, 38, 39, 42, 55, 56, 59, 63, 65, 69, \\ 70, 79, 81, 85, 89, 91, 130, 135, 137, 144, 146, 153, 156, 159, 161, 167, 170, 173, \\ 181, 190, 191, 192, 197, 
205, 207, 208, 241, 245, 246, 254, 262, 263, 266, 269
\}.$


\begin{table}[htb!]
\centering
\scalebox{.7}{
\begin{tabular}{|c|c|c|c|c|c|c|c|}
  \hline
No. & Arc & $|Aut(D)|$ & \# Par. cl. & \# Resol. &  $|Aut(C(D))|/|Aut(C(D^{\perp}))|$ \\
  \hline
1 & PG(2,16).1 & 68 & 2329 / 2329 & 409 / 409  & $2^2 17^1$ / $2^2 17^1 $ \\
  \hline
2 & PG(2,16).2 & 409 & 2550 / 2550 & 460 / 460 &  $2^3 3^1 17^1$ / $2^3 3^1 17^1$ \\
  \hline
3 & DEMP.1 & 24 & 250 / 319 & 52 / 52 & $ 2^{46} 3^{19} 5^9 7^6 11^3 13^3 17^3$ / $2^{33} 3^2$ \\
\hline
4 & DEMP.2 & 144 & 543 / 1023 & 52 / 214 & $2^{44} 3^{14}$ / $ 2^{18} 3^4 5^1 7^1 $  \\
  \hline
5 & DEMP.3 & 24 & 611 / 645 & 52 / 52 & $2^{45} 3^{13}$ / $ 2^{41} 3^{11} 5^2 7^2 $  \\
  \hline
6 & DEMP.4 & 48 & 531 / 691 & 52 / 52 & $2^{44} 3^{14}$ / $ 2^{45} 3^{15} $  \\
  \hline
7 & DEMP.5 & 4 & 255 / 377 & 52 / 52 & $ 2^{46} 3^{19} 5^9 7^6 11^3 13^3 17^3$ / $2^{33} 3^2$  \\
  \hline
8 & SEMI4.1 & 96 & 2569 / 2569 & 52 / 52 & $2^{17} 3^3$ / $ 2^{17} 3^3 $  \\
  \hline
9 & SEMI2.1 & 24 & 327 / 327 & 52 / 52 & $2^{45} 3^{15}$ / $2^{45} 3^{15}$  \\
  \hline
10 & SEMI2.2 & 144 & 1279 / 1279 & 55 / 55 & $2^{18} 3^4 5^1 7^1$ / $2^{18} 3^4 5^1 7^1$  \\
  \hline
11 & SEMI2.3 & 32 & 1497 / 1497 & 52 / 52 & $2^{26} 3^1$ / $2^{26} 3^1$  \\
  \hline
12 & SEMI2.4 & 32 & 1313 / 1313 & 52 / 52 &  $2^{25} 3^1 $ / $2^{25} 3^1$  \\
  \hline
13 & SEMI2.5 & 16 & 1045 / 1045 & 52 / 52 & $2^{37} 3^6 $ / $2^{37} 3^6$  \\
  \hline
14 & SEMI2.6 & 48 & 547 / 691 & 52 / 52 & $2^{45} 3^{15}$ / $2^{17} 3^3$  \\
  \hline
15 & SEMI2.7 & 48 & 691 / 547 & 52 / 52 & $2^{17} 3^3 $ / $2^{45} 3^{15}$  \\
  \hline
16 & LMRH.1 & 96 & 2265 / 2265 & 104 / 104 & $2^{45} 3^{15}$ /$2^{45} 3^{15}$  \\
 \hline
17 & LMRH.2 & 32 & 2377 / 2289 & 64 / 64 & $2^{45} 3^{15}$ /$2^{45} 3^{15}$  \\
 \hline
18 & MATH.1 & 24 & 291 / 275 & 52 / 52 &  $2^{49} 3^{20} 5^9 7^6 11^3 13^3$ / $2^{49} 3^{20} 5^9 7^6 11^3 13^3$ \\
\hline
19 & MATH.2 & 32 & 1729 / 1553 & 52 / 52 & $2^{37} 3^2$ / $2^{45} 3^{15}$  \\
\hline
20 & MATH.3 & 32 & 2401 / 2217 & 64 / 104 & $2^{45} 3^{15}$ / $2^{45} 3^{15}$ \\
  \hline
21 & MATH.4 & 32 & 1665 / 1473 & 52 / 52 & $2^{38} 3^5$ / $2^{38} 3^5$ \\
  \hline
22 & MATH.5 & 16 & 1233 / 1457 & 52 / 52 & $2^{43} 3^{12} 5^2 7^2$ / $2^{36} 3^4$ \\
  \hline
23 & MATH.6 & 16 & 1329 / 1405 & 52 / 52 & $2^{48} 3^{14} 5^6 7^6$ /$2^{45} 3^{15}$ \\
\hline
24 & MATH.7 & 16 & 1125 / 1505 & 52 / 52 & $2^{48} 3^{14} 5^6 7^6$ / $2^{45} 3^{15}$ \\
  \hline
25 & HALL.1 & 24 & 274 / 558 & 52 / 52 & $2^{49} 3^{20} 5^9 7^6 11^3 13^3$ / $2^6 3^2$ \\
  \hline
26 & HALL.2 & 4 & 309 / 445 & 52 / 52 & $2^{46} 3^{19} 5^9 7^6 11^3 13^3 17^3$ / $2^{15} 3^4$ \\
  \hline
27 & BBH1.1 & 24 & 330 / 330 & 52 / 52 & $2^{40} 3^{13} 5^3 7^3$ / $2^{40} 3^{13} 5^3 7^3$ \\
  \hline
28 & BBH1.2 & 32 & 2017 / 2017 & 136 / 136 & $2^{38} 3^5$ / $2^{38} 3^5$ \\
\hline
29 & BBH1.3 & 4 & 285 / 285 & 52 / 52 & $2^{43} 3^{16} 5^6 7^4 11$ / $2^{43} 3^{16} 5^6 7^4 11$ \\
\hline
30 & JOWK.1 & 16 & 1389 / 1241 & 52 / 52 & $2^{37} 3^2$ / $2^{44} 3^{14} 5^1 7^1$ \\
  \hline
31 & JOWK.2 & 32 & 2409 / 2321 & 104 / 52 & $2^{37} 3^2$ / $2^{38} 3^6 5^1$ \\
  \hline
32 & JOHN.1 & 32 & 1953 / 1641 & 144 / 52 & $2^{45} 3^{15}$ / $2^{48} 3^{14} 5^6 7^6$ \\
  \hline
33 & JOHN.2 & 32 & 1953 / 1841 & 144 / 52 & $2^{45} 3^{15}$ / $2^{48} 3^{14} 5^6 7^6$ \\
  \hline
34 & JOHN.3 & 32 & 2017 / 1761 & 136 / 52 & $2^{38} 3^5$ / $2^{48} 3^{14} 5^6 7^6$ \\
\hline
35 & JOHN.4 & 32 & 2409 / 1929 & 104 / 52 & $2^{37} 3^2$ / $2^{48} 3^{14} 5^6 7^6$ \\
  \hline
36 & DSFP.1 & 24 & 1045 / 1121 & 52 / 52 & $2^{45} 3^{13}$ / $2^{43} 3^{11}$ \\
  \hline
\end{tabular}}
\caption{Designs associated with maximal (52,4)-arcs}
\label{tab1}
\end{table}

\begin{table}[htb!] \centering
\scalebox{0.53}{
\begin{tabular}{|c|c|c|c|c|c|c|c|c|c|c|c|c|c|}
  \hline
   & PG(2,16) & DEMP & SEMI4 & SEMI2 & LMRH & MATH & HALL & BBH1 & JOWK & JOHN & DSFP & BBH2 & BBS4  \\
  \hline
PG(2,16) & & & & & & & {\bf 1} & & & & & & \\
DEMP & & & & {\bf 1,5} & {\bf 5} & {\bf 5} & {\bf 5} & & {\bf 3} & {\bf 5} & {\bf 2,5} & & \\
SEMI4 & & & {\bf 2} & {\bf 3,5} & {\bf 1} & & & & {\bf 1} & {\bf 1} & {\bf 1} & & {\bf 1} \\
SEMI2 & & {\bf 1,5} & {\bf 3,5} & & {\bf 5} & {\bf 5} & & & & {\bf 5} & & & \\
LMRH & & {\bf 5} & {\bf 1} & {\bf 5} & {\bf 4,5} & {\bf 5} & & & {\bf 2}  & {\bf 5} & {\bf 3} & & \\
MATH & & {\bf 5} & & {\bf 5} & {\bf 5} & {\bf 5} & {\bf 5} & & {\bf 2,5} & {\bf 5} & & & \\
HALL & {\bf 1} & {\bf 5} & & & & {\bf 5} & & {\bf 1} & & {\bf 1} & & {\bf 1} & \\
BBH1 & & & & & & & {\bf 1} & & & {\bf 4,5} & & & \\
JOWK & & {\bf 3} & {\bf 1} & & {\bf 2} & {\bf 2,5} & & & & {\bf 4,5} & & & \\
JOHN & & {\bf 5} & {\bf 1} & {\bf 5} & {\bf 5} & {\bf 5} & {\bf 1} & {\bf 4,5} & {\bf 4,5} & & & & \\
DSFP & & {\bf 2,5} & {\bf 1} & & {\bf 3} & & & & & & & & \\
BBH2 & & & & & & & {\bf 1} & & & & & & \\
BBS4 & & & {\bf 1} & & & & & & & & & & \\
  \hline
\end{tabular}}
\caption{Connections between projective planes of order 16}
\label{tab2}
\end{table}

\newpage

 \begin{table}[h!]
\centering
\scalebox{.7}{
\begin{tabular}{|r|c|l|l|c|}
  \hline
{\bf No.} & {\bf  2-rank} & {\bf (52,4)-Arc } & {\bf (A$_2$,A$_4$)} 
& $|Aut(C(D))|$\\
  \hline
1& 41 &  PG(2,16).1 & (0,221) &   $2^2 17^1$\\
\hline
2& 41 &  PG(2,16).2 & (0,221) & $2^3 3^1 17^1$  \\
  \hline
3& 43 & HALL.1$^\perp$ & (6,1037) & $2^6 3^2$ \\
  \hline
4& 45   &  $\{$DEMP.1$^\perp$,DEMP.5$^\perp$$\}$ & (24,3989) & $ 2^{33} 3^2$ \\
  \hline
5& 45 &  $\{$DEMP.2$^\perp$, SEMI2.2$\}$ & (6,4325) & $ 2^{18} 3^4 5^1 7^1 $ \\
\hline
6& 45   & $\{$SEMI4.1, SEMI2.7$\}$ & (0,4469) & $ 2^{17} 3^3 $ \\
\hline
7& 45   & SEMI2.3 & (18,4165) & $ 2^{26} 3^1 $ \\
\hline
8& 45   & SEMI2.4 & (16,4277) & $ 2^{25} 3^1 $ \\
\hline
9& 45   & HALL.2$^\perp$ & (12,4229) & $ 2^{15} 3^4 $ \\
\hline
10& 46   & JOHN.3 & (42,8293) & $2^{38} 3^5$ \\
\hline
11& 46   &  $\{ $JOHN.4, JOWK.1, MATH.2$\}$ & (26,8613) & $2^{37} 3^2$ \\
\hline
12& 46   & JOWK.2$^\perp$ & (46,8325) &  $2^{38} 3^6 5^1$ \\
  \hline
13& 46 & $\{$MATH.4, MATH.4$^\perp$$\}$ & (42,8549) & $2^{38} 3^5$ \\
\hline
14& 46   & MATH.5$^\perp$ & (42,8549) & $2^{36} 3^4$  \\
\hline
15& 46   & SEMI2.5 & (50,8453) & $2^{37} 3^6$  \\
\hline
16& 47   & BBH1.1 & (120,16853) &  $2^{40} 3^{13} 5^3 7^3$  \\
 \hline
17& 47   & $\{$DEMP.2, DEMP.4$\}$ & (72,17045) & $2^{44} 3^{14}$ \\
\hline
18& 47   & $\{$DSFP.1, DEMP.3$\}$ & (74,16997) & $2^{45} 3^{13}$ \\
\hline
19& 47   & DSFP.1$^\perp$ & (66,17093) &  $2^{43} 3^{11}$ \\
\hline
  &  & $\{$JOHN.1, LMRH.1, LMRH.2, LMRH.2$^\perp$, & &  \\
20& 47  & MATH.2$^\perp$, MATH.3, MATH.3$^\perp$, MATH.6$^\perp$, & (78,16901) & $2^{45} 3^{15}$ \\
&   &  MATH.7$^\perp$, SEMI2.1, SEMI2.6, DEMP.4$^\perp$$\}$ & & \\
\hline
21& 47   & JOWK.1$^\perp$ & (94,16709) &  $2^{44} 3^{14} 5^1 7^1$ \\
\hline
22 & 47   & MATH.5 & (106,16869) & $2^{43} 3^{12} 5^2 7^2$  \\
\hline
23 & 47   & DEMP.3$^\perp$ & (98,16965) & $2^{41} 3^{11} 5^2 7^2$  \\
\hline
24& 48 & $\{$JOHN.1$^\perp$,  JOHN.2$^\perp$, JOHN.3$^\perp$, & (174,33669) &  $2^{48} 3^{14} 5^6 7^6$  \\
&   &  JOHN.4$^\perp$, MATH.6, MATH.7$\}$ & &\\
\hline
25& 48 & BBH1.3 & (186,33829) & $2^{43} 3^{16} 5^6 7^4 11$  \\
\hline
26& 49   & $\{$HALL.1, MATH.1, MATH.1$^\perp \}$ & (366,67205) & $2^{49} 3^{20} 5^9 7^6 11^3 13^3$ \\
\hline
27& 49 & $\{$DEMP.1, DEMP.5, HALL.2$\}$ & (408,67541) & $ 2^{46} 3^{19} 5^9 7^6 11^3 13^3 17^3$  \\
  \hline
\end{tabular}}
\caption{Equivalence classes of codes}
\label{tab3}
\end{table}

  \begin{table}[h!]
\centering
\scalebox{.8}{
\begin{tabular}{|r|c|c|c|c|}
  \hline
{\bf No.} & {\bf Arc} & {\bf  [n,k,d] of $C(D)$} & {\bf [n,k$^\perp$,d$^\perp$] of $C(D)^\perp$ }  & $A_{d^\perp}$  \\
\hline
1& PG(2,16).1 & [52,41,4] & [52,11,18] & 54  \\
\hline
2& PG(2,16).2 & [52,41,4] & [52,11,18] & 54  \\
\hline
3& HALL.1$^\perp$ & [52,43,2] & [52,9,18] & 24 \\
\hline
4& DEMP.1$^\perp$ & [52,45,2] & [52,7,20] & 3  \\
\hline
5& DEMP.2$^\perp$ & [52,45,2] & [52,7,20] & 3 \\
\hline
6& HALL.2$^\perp$ & [52,45,2] & [52,7,20] & 3  \\
\hline
7& SEMI2.3 & [52,45,2] & [52,7,20] & 11  \\
\hline
8& SENI2.4 & [52,45,2] & [52,7,18] & 4  \\
  \hline
9& SEMI4.1 & [52,45,4] & [52,7,20] & 3  \\
\hline
10& JOHN.3 & [52,46,2] & [52,6,20] & 3  \\
  \hline
11& JOHN.4 & [52,46,2] & [52,6,20] & 3  \\
  \hline
12& JOWK.2$^\perp$ & [52,46,2] & [52,6,20] & 11  \\
  \hline
13& MATH.4 & [52,46,2] & [52,6,18] & 4 \\
\hline
14& MATH.5$^\perp$ & [52,46,2] & [52,6,18] & 4  \\
\hline
15& SEMI2.5 & [52,46,2] & [52,6,18] & 4  \\
\hline
16& BBH1.1 & [52,47,2] & [52,5,18] & 6  \\
  \hline
17& DEMP.2 & [52,47,2] & [52,5,20] & 3  \\
\hline
18& DEMP.3$^\perp$ & [52,47,2] & [52,5,18] & 4  \\
\hline
19& DSFP.1 & [52,47,2] & [52,5,20] & 3  \\
  \hline
20& DSFP.1$^\perp$ & [52,47,2] & [52,5,20] & 1  \\
  \hline
21& JOHN.1 & [52,47,2] & [52,5,20] & 3  \\
  \hline
22& JOWK.1$^\perp$ & [52,47,2] & [52,5,20] & 7  \\
  \hline
23& MATH.5 & [52,47,2] & [52,5,18] & 4  \\
\hline
24& BBH1.3 & [52,48,2] & [52,4,18] & 2  \\
  \hline
25& JOHN.1$^\perp$ & [52,48,2] & [52,4,20] & 3  \\
  \hline
26& DEMP.1 & [52,49,2] & [52,3,18] & 3  \\
\hline
27& HALL.1 & [52,49,2] & [52,3,20] & 3  \\
\hline
\end{tabular}}
\caption{Parameters of codes and their dual codes}
\label{tab4}
\end{table}

\begin{table}[ht!]
\centering
\caption{Equivalences between codes}
\label{tab5}
\scalebox{.8}{
\begin{tabular}{|c|l|}
  \hline
$(A,A')$ & $\pi: \ C(D(A)) \rightarrow C(D(A')))$  \\
\hline
&(1, 29, 13, 45)(2, 30, 16, 46, 4, 31, 14, 47)(3, 32, 15, 48) \\
(DEMP.1$^\perp$, DEMP.5$^\perp$) & (5, 24, 7, 22)(6, 21, 8, 23)(9, 36, 10, 33)(11, 35)(12, 34) \\
& (18, 19, 20)(25, 28, 27, 26)(37, 41)(38, 42)(39, 43)(40, 44)
 \\
\hline
&(1, 50, 37, 25, 17, 48, 8, 2, 51, 9, 11, 31, 13, 26,
 45, 12, 7, 14,\\
(SEMI2.2, DEMP.2$^\perp$) & 32, 39, 36, 18, 29, 40, 15, 10, 28, 35, 42, 6, 5)(3, 52, 46,
 22, \\ & 34, 27)(4, 49, 23, 24)(20, 21, 30)(38, 41, 43, 47) \\
\hline
& (2, 33, 17, 31, 24, 29, 49, 35, 46, 51)(3, 42, 32, 45, 22, 21, 13,\\
(SEMI4.1, SEMI2.7) & 26, 48, 14, 38, 9, 40, 50, 4, 11, 47, 44, 5, 30, 6, 41, 8, 23, 10,\\
& 18, 37, 25, 52)(7, 20, 39, 28, 34, 36, 12, 16, 27) \\
\hline
 & (1, 50)(2, 24, 21, 25, 26, 13, 9, 30, 35, 31, 40, 6, 51, 4, 42, 23,\\ 
(JOWK.1, MATH.2)& 43, 37, 5, 7, 49)(3, 20, 8, 10, 48, 36, 34, 14, 11, 28, 46, 52) \\ &
(12, 29, 33, 16, 38, 27, 18, 39, 45, 32, 19)(15, 17, 44)(22, 41)\\
  \hline
 & (1, 5, 7, 8, 12, 2, 6, 10)(3, 11)(4, 9)(13, 29, 25, 50, 41, 38, 14, \\
(JOHN.4, MATH.2) & 30, 28, 52, 37, 19, 32, 49, 16, 33, 26, 27, 24, 21, 47, 42, 18, 34,\\
& 45, 43)(15, 35, 46, 20, 36, 22, 48, 40, 17, 31, 51, 44, 39)\\
  \hline
 & (1, 40, 14, 11, 48, 24, 13, 7, 41, 31, 3, 49, 5, 47, 35, 36, 27, 45,\\
(MATH.4$^\perp$, MATH.4) & 33, 30, 4, 52, 19, 25, 44, 28, 38, 6, 39, 12, 51, 18, 21, 17, 29)\\
& (2, 46, 34, 32)(8, 37, 20, 26, 42, 23)(9, 50)(10, 43, 22)\\
  \hline
  & (2, 50)(3, 8, 51)(4, 12, 14, 52)(5, 28, 23, 20, 17, 21, 9, 7, 24, 15)\\
(DEMP.2, DEMP.4) & (6, 42, 31, 49, 22, 36, 29, 39, 16)(11, 19, 18, 47, 27, 48, 46, 43, \\
  & 44, 40, 35, 33, 13, 30, 32, 45, 37, 41, 38, 25, 34) \\
  \hline
  & (1, 49)(2, 51)(3, 24, 10, 6, 30, 25, 23, 21, 17, 22, 16, 41, 12, 26, \\
(DSFP.1,DEMP.3) & 44, 29, 40, 28, 32, 8, 13, 18, 52, 4, 7, 34, 11, 14, 43, 20, 33, 35, \\
  & 15, 50)(5, 48, 36, 27, 45, 39, 38, 31, 46, 42, 9, 47)(19, 37)\\
  \hline
  & (1, 9)(2, 6, 48, 47, 20, 42, 14, 36, 27, 16, 45, 17, 32, 40, 11, 3, \\
(LMRH.2, SEMI2.1) & 12, 4, 15, 23, 37, 35, 52, 34, 49, 38, 8, 21, 10)(5, 51, 44)(7, 18, \\
  & 39, 25)(13, 29, 19, 26)(22, 33)(24, 43, 28, 30)(31,
    46)\\
  \hline
  & (6, 8, 20, 45, 51, 14, 32, 28, 46, 12, 21, 30, 34, 38, 10, 9)(7, 17, \\
(LMRH.2$^\perp$, SEMI2.1)  & 36, 44, 16, 26, 40, 22, 33, 35, 41, 49, 50, 11, 18, 39, 19, 42, 52, \\
  & 47, 15, 23, 24, 27, 43, 13, 29, 31, 25, 37) \\
  \hline
  & (6, 11, 14, 28, 39, 7, 25, 18, 38, 13, 20, 44, 49, 27, 45, 46, 19, \\
(SEMI2.6, SEMI2.1) & 34, 23, 12, 47, 10, 35, 51, 22, 36, 26, 29, 21, 9, 17, 41, 43, 16,  \\
  & 50, 30, 15, 31, 48, 52, 33, 32, 42, 24)
\\
  \hline
  & (1, 25)(2, 28, 4, 31, 49, 12, 50, 15, 41, 37, 10, 14, 38, 7, 29, 13, 
\\
(DEMP.4$^\perp$, SEMI2.1) & 35, 17, 36, 20, 45, 9, 11, 47, 18, 39, 19, 42, 40, 22, 33, 5, 26)\\
  & (3, 34, 8, 32, 52, 51, 48, 21, 27)(6, 23, 24, 30, 16,
    44, 46)\\
  \hline
 & (6, 12, 17, 39, 44, 40, 32, 23, 45, 11, 24, 10, 27, 36, 29, 
16, 20,\\
(LMRH.1, SEMI2.1)& 7, 30, 35, 41, 37, 49, 9, 48, 47, 43, 21, 34, 38, 52, 18, 19, 25, 28,\\
& 22, 42, 14, 15, 33, 13, 51, 50, 46)(26, 31)\\
  \hline
 & (1, 29)(2, 32, 4, 23, 43, 41, 35, 28, 45, 5,
 6, 9, 12, 51, 19, 47, 17,\\
(MATH.2$^\perp$, SEMI2.1) & 11, 48, 20, 50, 10, 15, 24, 46, 8, 21, 37, 13, 30)(3, 26, 39, 49, 7,\\
& 18, 14, 27, 42, 38, 16, 33, 31)(22, 40, 52)(25, 36, 34) \\
  \hline
\end{tabular}}
\end{table}

\begin{table}[h!]
\ContinuedFloat  
\centering
\caption{Equivalences between codes (continued)}
\scalebox{.8}{
\begin{tabular}{|c|l|}
  \hline
$(A,A')$ & $\pi: \ C(D(A)) \rightarrow C(D(A')))$  \\
  \hline
 & (1, 5)(2, 8, 4, 11, 34, 21, 6)(3, 14, 20, 31, 32, 51, 52, 22, 15, 47,\\
(MATH.3, SEMI2.1) & 33, 42, 7)(9, 25)(10, 35, 45)(12, 41, 46, 40, 24, 36, 26, 23)(13, 50,\\
& 49, 27, 29, 48, 43, 19, 28, 39, 30, 18, 44)(16, 17, 38) \\
  \hline
 & (1, 29)(2, 32, 4, 23, 19, 41, 31, 3, 26, 39, 43, 28, 45, 5, 6, 9, 12,\\
(MATH.3$^\perp$, SEMI2.1)& 51, 49, 13, 30)(7, 18, 38, 40, 46, 8, 21)(10, 15, 24, 22)(11, 48, 20,\\
& 44, 34, 14, 27, 42, 25, 36, 50, 16, 33)(17, 35, 47)\\
  \hline
 & (1, 29)(2, 32, 4, 23, 41, 17, 7, 21, 27, 36, 14, 38, 24, 25, 39, 28,\\
(MATH.6$^\perp$, SEMI2.1) & 45, 43, 46, 47, 20, 37, 44, 50, 34, 8, 18, 5, 6, 9, 12, 48, 19, 11,\\
& 15, 13, 30)(3, 26, 42, 22, 52, 35, 40, 49, 16, 31)(10, 51, 33) \\
  \hline
 & (1, 29)(2, 32, 4, 23, 49, 34, 50, 33, 8, 18, 11, 51, 41, 46, 19, 5,\\
(MATH.7$^\perp$, SEMI2.1) & 6, 9, 12, 48, 20, 40, 52, 16, 38, 28, 45, 43, 14, 31, 3, 26, 42, 22,\\
& 25, 39, 24, 27, 36, 37, 44, 17, 7, 21, 35, 10, 15, 13, 30) \\
  \hline
& (1, 6)(2, 12, 48, 27, 43, 50, 30, 34, 36, 14, 35, 32, 42, 23, 40,\\
(JOHN.1, SEMI2.1) & 26, 41, 31, 11, 18, 13, 8, 4, 9, 51, 22, 52, 33, 25, 37, 28, 46, 7,\\
& 21, 17, 44, 45, 10)(5, 15, 38, 39, 47, 24, 20, 49, 19, 16) \\
  \hline
 & (1, 45)(2, 46)(3, 47)(4, 48)(5, 29, 42, 40, 25, 24, 36, 37, 17, 18,\\
(JOHN.1$^\perp$, MATH.6) & 26, 13, 9)(6, 30, 39, 20, 22, 34, 43, 49, 27, 14, 10)(7, 31, 50,\\
& 28, 16, 12, 8, 32, 51, 23, 35, 44, 52, 15, 11)(19, 21, 33, 41, 38) \\
\hline
 & (1, 45)(2, 46)(3, 47)(4, 48)(5, 29, 42, 40, 25, 24, 36, 37, 17, 18,\\
(JOHN.2$^\perp$, MATH.6)& 26, 13, 9)(6, 30, 39, 20, 22, 34, 43, 49, 27, 14, 10)(7, 31, 50,\\
& 28, 16, 12, 8, 32, 51, 23, 35, 44, 52, 15, 11)(19, 21, 33, 41, 38) \\
\hline
& (1, 37)(2, 38)(3, 40, 4, 41, 33, 29, 15, 24, 52, 8, 46, 10, 17, 47,\\
(JOHN.3$^\perp$, MATH.6) & 11, 19, 50, 6, 42, 34, 30, 25, 13, 21, 43, 35, 31, 27, 16, 18,\\
& 48, 12, 20, 51, 7, 45, 9, 23, 49, 5, 39)(14, 22, 44, 36, 32, 28, 26) \\
  \hline
 & (1, 37)(2, 38)(3, 40, 4, 41, 33, 29, 15, 24, 52, 8, 46, 10, 17, 47,\\
(JOHN.4$^\perp$, MATH.6) & 11, 19, 50, 6, 42, 34, 30, 25, 13, 21, 43, 35, 31, 27, 16, 18,\\
& 48, 12, 20, 51, 7, 45, 9, 23, 49, 5, 39)(14, 22, 44, 36, 32, 28, 26) \\
  \hline
  & (6, 10, 13, 16, 21, 17)(7, 19, 23, 28, 22, 11, 36, 48, 49, 20, 47, \\
(MATH.7, MATH.6) & 35, 43, 40)(8, 29, 39, 32, 37, 26)(9, 27, 
    34, 18, 45, 31, 30,
\\
  & 24, 14, 42)(12, 33, 46, 44, 15, 38, 51, 52, 25, 50, 41)\\
  \hline
 & (1, 33)(2, 7, 47, 38, 52, 41, 32, 24, 13, 36, 4, 25, 19, 9, 34)\\
(MATH.1$^\perp$, MATH.1) & (3, 8, 12, 28, 15, 26, 20, 17, 46, 37, 49, 16, 10, 35)(5, 27, 50\\
&  31, 22, 23, 51, 39, 14, 6, 45, 30, 18, 48, 43, 42, 40, 29, 11) \\
  \hline
& (6, 10, 13, 16, 21, 17)(7, 19, 23, 28, 22, 11, 36, 48, 49, 20, 47,\\
(HALL.1, MATH.1) & 35, 43, 40)(8, 29, 39, 32, 37, 26)(9, 27, 34, 18, 45, 31, 30,\\
& 24, 14, 42)(12, 33, 46, 44, 15, 38, 51, 52, 25, 50, 41) \\
\hline
  & (1, 51, 37, 46, 31, 36, 43, 28, 27, 16, 26, 52, 50, 49, 41, 23, 25, \\
(DEMP.1, HALL.2) & 47, 32, 40, 4, 22, 3, 2)(5, 7, 13, 24)(6, 33, 10, 42, 18, 17, 35, 
\\
  & 12, 9, 8, 19, 39, 14, 11, 15, 29)(20, 45, 21)(38, 48,
    44)\\
  \hline
  & (1, 49, 51, 48, 16, 34, 40, 37, 4, 19, 46, 42, 7, 26, 25, 6, 2, 17, \\
(DEMP.5, HALL.2) & 23, 11, 5, 24, 22, 39, 32, 45, 41, 50, 52, 47, 15, 31, 36, 29, 44, \\
  & 38, 8, 3)(9, 28, 27, 10, 30, 33, 18, 35)(20, 43)
\\
  \hline

\end{tabular}}
\end{table}

\end{document}